\documentclass[11pt,fleqn,twoside]{article}
\usepackage{amsfonts,amssymb,latexsym}
\makeatletter
\newcommand{\prava}{\footnotesize\it
\begin{flushright}
\begin{minipage}{18cm}
Copyright \copyright 1998 by A.A. Borghardt and D.Ya. Karpenko
\end{minipage}
\end{flushright}}

\newcommand{\name}[1]{\begin{flushleft}
                       \LARGE \bf #1
                       \end{flushleft}\vspace{-3mm}}

\newcommand{\Author}[1]{\begin{flushleft}
                       \it #1 \end{flushleft}}

\newcommand{\Adress}[1]{\begin{flushleft}
                       \it #1 \end{flushleft}}

\newcommand{\Date}[1]{\begin{flushleft}
                      \small  \it #1 \end{flushleft}}

\newcommand{\ehkol}{Author \ name}
\newcommand{\ohkol}{Article \ name}
\renewcommand{\@evenhead}{
\hspace*{-3pt}\raisebox{-15pt}[\headheight][0pt]{\vbox{\hbox to \textwidth
{\thepage \hfil \ehkol}\vskip4pt \hrule}}}
\renewcommand{\@oddhead}{
\hspace*{-3pt}\raisebox{-15pt}[\headheight][0pt]{\vbox{\hbox to \textwidth
{\ohkol \hfil \thepage}\vskip4pt\hrule}}}
\renewcommand{\@evenfoot}{}
\renewcommand{\@oddfoot}{}

\newcommand{\be}{\begin{equation}}
\newcommand{\ee}{\end{equation}}
\newcommand{\ba}{\hspace*{-5pt}\begin{array}}
\newcommand{\ea}{\end{array}}
\newcommand{\p}{\partial}
\newcommand{\ds}{\displaystyle}
\makeatother

\font\BoldMath=cmmib10 scaled \magstep1

\newcommand{\pmA}{\mbox{\BoldMath \char 65}}
\newcommand{\pmk}{\mbox{\BoldMath \char 107}}
\newcommand{\pmq}{\mbox{\BoldMath \char 113}}

\begin{document}
\setcounter{page}{357}
\thispagestyle{empty}

\renewcommand{\ehkol}{A.A. Borghardt and D.Ya. Karpenko}
\renewcommand{\ohkol}{Fundamental  Solution of the Volkov Problem 
(Characteristic Representation)}

\begin{flushleft}
\footnotesize \sf
Journal of Nonlinear Mathematical Physics \qquad 1998, V.5, N~4,
\pageref{borghart-fp}--\pageref{borghart-lp}.
\hfill {\sc Letter}
\end{flushleft}

\vspace{-5mm}

\renewcommand{\footnoterule}{}
{\renewcommand{\thefootnote}{} \footnote{\prava}} 

\name{Fundamental  Solution of the Volkov Problem 
(Characteristic Representation)}\label{borghart-fp}

\Author{A.A. BORGHARDT and  D.Ya. KARPENKO}

\Adress{Donetsk Physical\/-Technical Institute, Ukrainian National
Academy of Sciences,  \\
340114 Donetsk, Ukraine \\
E-mail: borg@host.dipt.donetsk.ua}

\Date{Received June 10, 1998; Accepted July 13, 1998}

\begin{abstract}
\noindent 
The characteristic representation, or Goursat problem, for the
Klein-Fock-Gordon equation with Volkov interaction [1] is
regarded. It is shown that in this representation the explicit form
of the Volkov propagator can be obtained. Using the characteristic
representation technique, the
Schwinger integral [2] in the Volkov problem can be calculated. 
\end{abstract}

\section{Introduction}

The Klein-Fock-Gordon (KFG) equation is the basic equation of
relativistic quantum mechanics (RQM).
From the Cauchy problem point of view this equation  has the 
advanced or retarded Green function as well as the causal or Feynman 
Green function $\Delta_C$. In this paper our attention is focused on 
the  function $\Delta_C$, i.e., 
\[
\Delta_C=\frac 12 \left(\Delta_S +i \Delta^{(1)}\right),
\]
where $\ds \frac 12 \Delta_S$ is the real part of the Feynman propagator 
$\Delta_C$, and is written in the form
\be 
\Delta_S =\frac{1}{2\pi} \left( \delta\left(\lambda^2\right) -
\theta \left(\lambda^2\right) k_0 \frac{J_1(k_0\lambda)}{\lambda}\right).
\ee
The following notation is used:  $\lambda =\sqrt{c^2t^2 -r^2} =
\sqrt{x_\mu^2} $ is the 4-D interval of Minkowski space $M_{(+)}^4$, 
$\theta$ is the Heaviside step function, $k_0 =\mu c/\hbar$ 
is the inverse Compton wave length, and $J_1$ is the Bessel function of 
the f\/irst kind with index 1.

	The function $\Delta^{(1)}$ is the imaginary part of 
the  function $\Delta_C$, i.e.,
\be
\Delta ^{(1)}= \frac{k_0}{4\pi} \left\{
\theta \left(\lambda^2\right) \frac{N_1(k_0\lambda)}{\lambda} +
\frac{2}{\pi} \theta(\widetilde \lambda^2) 
\frac{K_1(k_0\widetilde \lambda)}{\widetilde \lambda}
\right\},
\ee
$\widetilde \lambda =\sqrt{r^2-c^2 t^2} =\sqrt{-x_\mu^2}$ 
is the 4-D interval of Minkowski space $M_{(-)}^4$ $(\widetilde
\lambda\geq 0)$, $N_1$ 
is the Bessel function of the second kind, and $K_1$ 
is the MacDonald function.

In free case from, the  function $\Delta_C$  can be obtained 
for all propagators RQM.

It should be noted that the explicit form of the  function $\Delta_C$
is known for the free KFG equation only. The causal propagator $\Delta_C$
has not been calculated, even for the simplest interaction which is the Volkov 
interaction.

The purpose of this letter is to calculate the Green function for the 
Volkov interaction from the CR viewpoint and to seek 
the  function $\Delta_C$ with the Volkov interaction from 
Cauchy viewpoint.

\section{Free KFG equation from CR viewpoint}

\setcounter{equation}{0}

In the CR viewpoint we require an additional spatial coordinate.
For the sake of def\/inition, let  $z$ be this coordinate. 
The wave equation  may then be rewritten in the form 
\be
\left( \frac{\p^2}{\p z^2} -\frac{1}{c^2} \frac{\p^2}{\p t^2} -
\hat a^2\right) \Psi (x_1, x_2, z,t)=0,
\ee
where the operator parameter $\hat a^2$ does not contain  second 
derivatives with respect to $t$ and $z$.

For the KFG equation the operator parameter $\hat a^2$
is given by
\be 
\hat a^2 =k_0^2 -\Delta_\perp, \qquad 
\Delta_\perp = \frac{\p^2}{\p x_1^2} +\frac{\p^2}{\p x_2^2}.
\ee
If  $\hat a^2$ is a positive constant we obtain the 1-D 
telegraph equation
\be 
\left( \hat L -a^2\right) \Psi(z,t)=0,
\ee
where    $\ds \hat L =\frac{\p^2}{\p z^2} -\frac{1}{c^2} \frac{\p^2}{\p t^2}=
-4 \frac{\p^2}{\p \xi \p \eta}$ is the 1-D wave operator and  
$\xi =ct-z$, $\eta=ct +z$ are the characteristics of a wave equation.

The fundamental solution for equation (2.3) is well known. 
It is called the Riemann function~[4]
\be
\Psi (z,t) =J_0 \left( \sqrt{\left( c^2 t^2 -z^2\right) a^2}\right)=
J_0\left( \sqrt{\xi \eta a^2}\right),
\ee
where $J_0$  is the Bessel function.

We are interested in the CR for the KFG equation (see Refs. [3, 4] 
and references therein). In this case equation (2.4) has the form 
\be 
\Psi(x_1, x_2,z,t) =J_0\left( \sqrt{\left( c^2 t^2-z^2\right) \hat a^2}
\right) |0\rangle,
\ee
where we take the operator parameter  from equation (2.2).
Here $|0\rangle$ is the initial value of the wave function along 
the characteristics $|t| = |z|/c$ (or $ct\pm z = 0$).

The solution (2.5) is stable, provided that the value under the square 
root is positive. 
Otherwise, if the value under the square 
root is negativ, the solution will be increasing.

Equation (2.5), with the operator  $\hat a^2$, should be considered  
in two regions of the Min\-kow\-s\-ki space $M^2$. The space $M^2$ consist of 
the space $M^2_{(+)}$, where $c^2t^2 - z^2\geq 0$ (or $|t| \geq |z|/c$),
and the space $M^2_{(-)}$, where $c^2t^2 - z^2\leq  0$ (or $|t|
\leq |z|/c$). 

We f\/irst consider the region $c^2t^2 - z^2 \geq  0$. As the initial 
value $|0\rangle$ we select the 2-D $\delta$-function, 
namely when $|t| = |z|/c$.  We have the initial localized function~[3]
\be 
|0\rangle =\frac{1}{4\pi^2} \int \exp (ik_1 x_1 +i k_2 x_2) \, dk_1 \,dk_2=
\delta(x_1,x_2).
\ee
Substituting equation (2.6) into equation (2.5), and integrating 
over the polar angle with the help of the formula 
\[
\int\limits_0^{2\pi} \exp(ik x\cos (\varphi-a)) d\varphi= 2\pi J_0(kx),
\]
 we get the fundamental solution  in the CR in the   space 
$M^2_{(+)}$. Therefore, we denote the wave function of equation (2.5) 
as $\Psi^{(+)}$
\be
\Psi^{(+)} (x_1,x_2,z,t) =\frac{1}{2\pi} \int\limits_0^\infty
J_0 \left( \sqrt{\left(c^2 t^2 -z^2\right) \left( k^2+k_0^2\right)}\right)
J_0(k x_\perp) k \, dk,
\ee
where $x_\perp =\sqrt{x_1^2 +x_2^2}$.

It is easy to calculate the integral (2.7) if $k_0\to 0$. 
In this case we have the photon propagator~[3]
\[
\Psi^{(+)}\Bigl|_{k_0\to 0} =\frac{\delta \left(\sqrt{c^2t^2 -z^2} -x_\perp
\right)}{2\pi x_\perp} =\frac{\delta(|t|-r/c)}{2\pi cr}=
\frac{1}{\pi} \delta\left(\lambda^2\right).
\]

	For $k_0 \not= 0$, we calculate the equation (2.7) 
with the help of the Sonin discontinuous integral
\be
\int\limits_0^\infty k^{n+1} \frac{J_m \left(\tau \sqrt{k^2+k_0^2}\right)}
{\left( k^2+k_0^2\right)^{m/2}} J_n(kx_\perp) \, dk =
\theta (\tau-x_\perp) \frac{x_\perp^n}{\tau^m}
\left( \frac{ \lambda}{x_\perp}\right)^{m-n-1} \!\!\!\!\!
J_{m-n-1}(k_0 \lambda),
\ee
where $\lambda =\sqrt{\tau^2 -x^2_\perp}$, and  $m > n > 0$.

Note that the integral (2.8) cannot be used directly,
because in equation (2.7) 
numbers  $m = n = 0$. It is therefore necessary to increase the order of the 
Bessel function. In order to achieve this we use the  equality from 
the Bessel function theory
\be 
J_0 \left( \tau \sqrt{ k^2+k_0^2}\right)=
\frac{1}{\tau} \frac{d}{d\tau} \left( \tau
\frac{J_1\left( \tau \sqrt{k^2+k_0^2}\right)}{\sqrt{k^2+k_0^2}} \right).
\ee

Substituting equation (2.9) into equation (2.7),  and taking 
into account equation (2.8) at $m=1$ and $n=0$, we obtain the
fundamental solution
of 
the KFG equation in the space $M^2_{(+)}$ as
\be
\ba{l}
\ds \Psi^{(+)} =\frac{1}{2\pi} \left( \frac{\delta(\tau-x_\perp)}{x_\perp}-
k_0 \theta (\tau -x_\perp) \frac{J_1(k_0 \lambda)}{\lambda}\right)
\vspace{3mm}\\
\ds \qquad =
\frac{1}{\pi} \left( \delta\left(\lambda^2\right) -
\frac{k_0}{2} \theta\left(\lambda^2\right)
\frac{J_1(k_0\lambda)}{\lambda}\right), 
\ea
\ee
where  $\tau =\sqrt{c^2t^2-z^2}=\sqrt{\xi \eta}$ 
is the timelike interval of the space $M^2_{(+)}$.

	The fundamental solution (2.10) from the CR viewpoint coincides with the function 
$\Delta_S$ of equation (1.1) from the Cauchy problem viewpoint.

We should also consider the fundamental solution in the space $M^2_{(-)}$,   where 
$z^2- c^2t^2 \geq  0$  (or $|t|\leq |z|/c$). 
For this purpose we rewrite equation~(2.5)~[4] as
\be 
\Psi(x_1,x_2,z,t) =J_0 \left( \sqrt{\left( z^2-c^2t^2\right) \left(
\Delta_\perp -k_0^2\right)}\right) |0\rangle.
\ee
We denote the wave function of equation (2.11) in the space $M^2_{(-)}$
 as  $\Psi^{(-)}$. It should be noted that the initial value
$|0\rangle$,   on the characteristics $|t| = |z|/c$,
cannot be selected in the localized form since we obtain 
unstable solutions~[3]. Therefore we take the superposition of the MacDonald 
functions  $K_0$~[3]
\be
|0\rangle =\int\limits_{k_0}^\infty K_0 (k x_\perp) k \, dk=k_0
\frac{K_1 (k_0 x_\perp)}{2\pi x_\perp} \ {\mathop {\longrightarrow}\limits_{
\mbox{\scriptsize $k_0\to 0$}}}
\ \frac{1}{2\pi x_\perp^2}.
\ee

Taking into account the equality on eigenfunctions $\Delta_\perp K_0=k^2K_0$,
and substituting equation (2.12) into equation (2.11), we obtain the 
fundamental solutions of the KFG equation in the space~$M^2_{(-)}$~[3]
as
\be 
\ba{l}
\ds \Psi^{(-)} (x_1.x_2,z,t) =\frac{1}{2\pi} \int\limits_{k_0}^\infty
J_0\left( \sqrt{ \left(z^2-c^2t^2\right) \left(k^2-k_0^2\right)}\right)
K_0(k x_\perp)k\, dk
\vspace{3mm}\\
\ds \qquad = \frac{k_0}{2\pi} \frac{K_1(k_0\widetilde \lambda)}{\widetilde
\lambda} \ {\mathop {\longrightarrow }\limits_{\mbox{\scriptsize $k_0\to 0$}} } \ \frac{1}{2\pi}
\frac{1}{\widetilde \lambda^2} >0.
\ea
\ee
The function $\Psi^{(-)}$ of equation (2.13), in the CR viewpoint, does not
coincide with the $\Delta^{(1)}$  function (see equation (1.2)).

\section{Volkov problem from CR viewpoint}

\setcounter{equation}{0}

As an example of the CR application, we regard the KFG equation with
electromagnetic interaction 
\be
\left(\hat L+\Delta_\perp -k_0^2 +\frac{2ie}{\hbar c} A_\mu \p_\mu -
\frac{e^2}{\hbar^2 c^2} A_\mu^2\right)\Psi =0,
\ee
where Volkov potential has the form 
\be
A_\mu =(A_1(\xi), A_2(\xi), 0,0), \qquad \xi=ct-z.
\ee
The wave function of equation (3.1), with the potential (3.2), can be expanded
into a Fourier integral over transverse coordinates $x_1$, $x_2$, i.e., 
\be
\Psi(x_1,x_2,z,t) =\frac{1}{4\pi^2} \int \exp(ik_1x_1 +ik_2 x_2)
\Phi (k_1,k_2,z,t) \, dk_1 \, dk_2.
\ee
Taking into account the transformation (3.3), we obtain the 1-D telegraph
equation with the variable parameter $K^2(\xi)$
\be
\left( 4\frac{\p^2}{\p \xi \p \eta} +K^2(\xi)\right)\Phi(k_1,
k_2,\xi,\eta)=0, 
\ee 
where 
\be
K^2(\xi) =k^2+k_0^2+\frac{2e}{\hbar c} ({\pmA \pmk}) +\frac{e^2}{\hbar^2
c^2} A_\mu^2, \qquad k^2=k_1^2+k_2^2.
\ee

Since the characteristic variable (in equation (3.4) is cyclic, we can seek
the solution of equation (3.4) with the initial condition on the
characteristics 
\be
\Phi(k_1,k_2,0,0)=1
\ee
in the form 
\be 
\Phi =J_0 \left( \sqrt{\eta f(\xi)}\right),
\ee
where the unknown function $f(\xi)$ satisf\/ies to the equation $f'(\xi) =
K^2(\xi)$. Whence it follows that 
\be
f(\xi) =\int\limits_0^\xi K^2(z) dz.
\ee

To prove equation (3.8) we need the following property $\hat L$-operator 
\be
4\frac{\p^2}{\p \xi \p \eta } J_0 \left( \sqrt {\eta f(\xi)}\right) =
-f' (\xi) J_0\left( \sqrt{\eta f(\xi)}\right).
\ee

We still need to calculate the integral (3.3) with the function $\Phi$,
which is determined by equations (3.7) and (3.8), i.e.,
\be
\ba{l}
\ds \Psi =\frac{1}{4\pi^2} \int dk_1\, dk_2 \exp i (k_1 x_1 +k_2x_2)
J_0 \left( \sqrt{ \xi \eta \langle K^2(\xi)\rangle }\right) ,
\qquad \mbox{and}
\vspace{3mm}\\
\ds \langle K^2(\xi)\rangle =\frac{1}{\xi} \int\limits_0^\xi K^2(z) \, dz.
\ea
\ee

The integral (3.10) reduces to the free case if we make the shift of the
momentum 
\[
{\pmq }={\pmk} +\frac{e}{\hbar c}\langle {\pmA}\rangle.
\]
In fact, the 2-D vector ${\pmq}$ is the general momentum. After this we
get the same integral as in the free case but the value $k_0$ is replaced by
$k_0(\xi)$
\be
\Psi =\!\frac{1}{4\pi^2} \exp\! \left( -\frac{ie}{\hbar c} x_\mu \langle A_\mu
\rangle \right)\!\! \int \exp i (q_1 x_1 +q_2 x_2) J_0\!
\left(\! \sqrt{\xi \eta \left( q^2+k_0^2(\xi)\right)}\!\right)\! dq_1  dq_2,
\ee
where
\be
k_0(\xi) =k_0 \sqrt{1+\left(\frac{e}{\hbar ck_0}\right)^2
\left(\langle A_\mu^2\rangle -\langle A_\mu\rangle^2\right)}.
\ee
Reasoning as in the calculation of the integral (2.7), we obtain the
fundamental solution of
the KFG equation with the Volkov interaction from the CR viewpoint 
\be
\Psi =\frac{1}{\pi} \exp \left(-\frac{ie}{\hbar c} x_\mu \langle
A_\mu\rangle \right) \left\{ \delta\left(\lambda^2\right) -
\frac{k_0(\xi)}{2} \theta \left(\lambda^2\right) 
\frac{J_1(k_0(\xi)\lambda)}{\lambda}
\right\}.
\ee

The solution (3.13) allows us to conclude that the particle in the
electromagnetic wave becomes more heavy. Its mass is increased by 
$\ds \sqrt{1+\left( \frac{e}{\hbar c k_0}\right)^2 \left( \langle
A_\mu^2\rangle -\langle A_\mu\rangle^2\right)^2}$ times. The
same conclusion follows from [5], where the Volkov problem is considered from
the classical Lorentz equation, and also from [6] (page 188; The nonlinear
interaction of massless spinor particle with electromagnetic f\/ield generates
the mass of a particle). 

The propagator of the Dirac equation with the Volkov interaction was considered
in  papers [2, 7, 8] as well as in many other papers. In
particular in [2] the propagator was found in the form of the proper time
integral. If we neglect the spin dependence and perform the change of
variable in integration $\alpha = 1/4s$, then the propagator determined by
equation~(4.26) in [2] can be represented in the form 
\be
G_s (x',x'') =-\frac{1}{4\pi^2} C(x',x'') \int\limits_0^\infty \exp
\left( \frac{-ik_0^2 (\xi',\xi'')}{4\alpha} -i\lambda \alpha \right) d\alpha,
\qquad \lambda^2 =-x_\mu^2,
\ee
where  
\[
k_0^2 (\xi', \xi'') =k_0^2 +\frac{1}{\xi'-\xi''} \int\limits_{\xi''}^{\xi'}
A_\mu^2 (\xi) \, d\xi -\left( \frac{1}{\xi'-\xi''}
\int\limits_{\xi''}^{\xi'} A_\mu (\xi) \, d\xi\right)^2.
\]
The integral (3.14) was calculated in the appendix of Schwingers paper [9]. 
The real part of the integral (3.14) is equal to the expression in
brackets
$\{\; \} $ of equation~(3.13). The imaginary part of the integral (3.14) is
equal to the function $\Delta^{(1)}$ of equation~(1.2), in which $k_0 \to
k_0(\xi )$. 

It is important to note that the factor $C(x',x'')$ coincides with the
factor before the integral (3.11). The latter was found by us while
considering the Volkov problem from the CR viewpoint. 

Thus, up to the phase factor, the Schwinger propagator (3.14) for the Volkov
interaction~(3.2), is the free function $\Delta_C$ with $k_0(\xi)$ 
\be
G_S (x',x'') =-\frac{1}{4\pi^2} \exp\!\left( \!-i(x'-x'')_\mu
\frac{1}{\xi'-\xi''} \!\int\limits_{\xi''}^{\xi'} \!\!A_\mu (\xi)
 \, d\xi\right) \! 
\Delta _C(x',x'',k_0(\xi',\xi'')).
\ee

 \label{borghart-lp}
\end{document}